\documentclass[12pt]{article}
\usepackage{amsmath}
\usepackage{amsfonts}
\usepackage{amssymb}
\usepackage{amsthm}

\usepackage{hyperref}

\newcommand{\Rational}{\mathbb{Q}}
\newcommand{\Irrational}{\Real\backslash\Rational}
\newcommand{\Natural}{\mathbb{N}}
\newcommand{\Neotr}{\Natural_0}
\newcommand{\UD}{\mathbb{D}}
\newcommand{\UC}{\mathbb{T}}
\newcommand{\Complex}{\mathbb{C}}
\newcommand{\ComplexE}{\overline{\mathbb{C}}}
\newcommand{\Real}{\mathbb{R}}

\newcommand{\mIm}{\Im {\rm m}\,}
\newcommand{\mRe}{\Re {\rm e}\,}

\newcommand{\Julia}{\mathcal J}
\newcommand{\Fatou}{\mathcal F}
\newcommand{\Basr}{\mathcal{A}^*}
\newcommand{\Linija}{\mathcal L}
\newcommand{\krug}[2]{\mathrm{D}(#1,#2)}
\newcommand{\res}{\mathop{\mathrm{R\,e\,s}}}
\newcommand{\lfr}{\hbox{\mathversion{bold}$\{$}}
\newcommand{\rfr}{\hbox{\mathversion{bold}$\}$}}

\newcommand{\dist}{\mathrm{dist}}

\newtheorem{theorem}{Theorem}
\newtheorem{result}{Theorem}
\newtheorem{lemma}{Lemma}
\newtheorem{proposition}{Proposition}

\theoremstyle{remark}
\newtheorem{remark}{Remark}
\newtheorem{example}{Example}

\newcommand{\mcite}[1]{\csname b@#1\endcsname}

\begin{document}
\title{Carath\'eodory convergence of immediate basins of attraction to a Siegel disk}
\author{Pavel Gumenuk${}^\dagger$}
\date{}

\maketitle {\let\thefootnote=\relax\footnotetext{\hskip-1.6em 2000 {\it Mathematics Suject Classification}.
Primary 30D05, 37F45. Secondary 37F50.\\[.3ex]{\it Key words and phrases.} Iteration of
analytic functions, Fatou set, Siegel disk, basin of attraction, convergence as to the
kernel.\\[.5ex]${}^\dagger$This work is partially supported by the {\it Research Council of Norway}, the {\it Russian Foundation for Basic Research} (grant \#07-01-00120),  and {\it ESF Networking Programme} ''Harmonic and Complex Analysis and its Applications''}

\maketitle

\begin{abstract}Let $f_n$ be a sequence of analytic functions in a domain $U$ with a common
attracting fixed point~$z_0$. Suppose that $f_n$ converges to $f_0$ uniformly on each compact
subset of $U$ and that $z_0$ is a Siegel point of $f_0$. We establish a sufficient condition
for the immediate basins of attraction $\mathcal A^*(z_0,f_n,U)$ to form a sequence that
converges to the Siegel disk of $f_0$ as to the kernel with respect to $z_0$. The same
condition is shown to imply the convergence of the K\oe nigs functions associated with~$f_n$ to
that of~$f_0$. Our method allows us also to obtain a kind of quantitative result for analytic one-parametric families.
\end{abstract}

\section{Introduction}

\subsection{Preliminaries} Let $U$ be a domain on the Riemann sphere $\ComplexE$ and
$f:U\to\ComplexE$ a meromorphic function. Define~$f^n$, the $n$-fold {\it iterate} of $f$, by
the following relations: ${f^1: U\to\ComplexE}$, ${f^1:=f}$,
${f^{n+1}:\big(f^n\big)^{-1}(U)\to\ComplexE}$, ${f^{n+1}:=f\circ f^n}$, ${n\in\mathbb N}$. It
is convenient to define $f^0$ as the identity map of~$U$. Denote%
\[E(f,U):=\bigcap_{n\in\Natural}\big(f^n\big)^{-1}(U).\]
The {\it Fatou set} $\Fatou(f,U)$ of the function $f$ (with respect to the domain~$U$) is the
set of all interior points $z$ of $E(f,U)$ such that  $\{f^n\}_{n\in\Natural}$ is a normal
family in some neighbourhood of $z$. Define the {\it Julia set} $\Julia(f,U)$ of $f$ (with
respect to the domain~$U$) to be the complement $U\setminus\Fatou(f,U)$ of the Fatou set.

Classically iteration of analytic (meromorphic) functions has been studied for the case of
$U\in\big\{\ComplexE,\Complex,\Complex^*:=\Complex\setminus\{0\}\big\}$ and $f:U\to U$, see
survey papers~\cite{AnIntro,ErLyub} for the details. As an extension the cases of
transcendental meromorphic functions and functions meromorphic in $\ComplexE$ except for a
compact totally disconnected set have been also investigated, see
e.g.\,\cite{Meromorph,Baker}. (Note that $f(U)\not\subset U$ for these cases.) In this paper
we shall restrict ourselves by the following

\noindent{\bf Assumption.} Suppose that $U, f(U)\subset\Complex$, i.e. $f$ is an analytic
function in a subdomain~$U$ of $\Complex$.

One of the basic problems in iteration theory of analytic functions is to study how the limit
behaviour of iterates changes as the function $f$ is perturbed. A large part of papers in this
direction are devoted to the continuity property  for the dependence of the Fatou and Julia
sets on the function to be iterated. We mention the work of A.\,Douady~\cite{Douady}, who
investigates the mapping $f\mapsto \Julia(f,\Complex)$ from the class of polynomials of fixed
degree to the set of nonempty plane compacta equipped with the Hausdorff metric
$d_H(X,Y):=\max\{\partial(X,Y),\partial(Y,X)\}$, $\partial(X,Y):=\sup_{x\in X}\dist(x,Y).$ We
also mention subsequent papers~\cite{Kisaka,Kriete1,Kriete2,Kriete3,Kriete4,Wu}
dealing with other classes of functions. Continuity of Julia sets is closely related to
behaviour of connected components of the Fatou set containing periodic points. Now we recall
necessary definitions.

Let $z_0\in U$ be a fixed point of $f$.  The number $\lambda:=f'(z_0)$ is called the {\it
multiplier} of $z_0$. According to the value of $\lambda$ the fixed point~$z_0$ is  said to be
{\it attracting} if $|\lambda|<1$, {\it neutral} if $|\lambda|=1$, and {\it repelling} if
$|\lambda|>1$. An attracting fixed point is {\it superattracting} if $\lambda=0$, or {\it
geometrically attracting} otherwise. Suppose $z_0$ is a neutral fixed point of~$f$ and none of
$f^n$, $n\in\Natural$, turns into the identity map; then the fixed point~$z_0$ is {\it
parabolic} if $\lambda=e^{2\pi i\alpha}$ for some $\alpha\in\Rational$, or {\it irrationally
neutral} otherwise. If an irrationally neutral fixed point belongs to $\Fatou(f,U)$, then it
is called a {\it Siegel point}.

The component of the Fatou set~$\Fatou(f,U)$ that contains a fixed point~$z_0$ is called the
{\it immediate basin} of $z_0$ and denoted by~$\Basr(z_0,f,U)$. The immediate basin of a
Siegel point is called a {\it Siegel disk}, and the immediate basin of an attracting fixed
point is called an {\it immediate basin of attraction}. It is a reasonable convention to put
by definition $\Basr(z_0,f,U):=\{z_0\}$ for fixed points~$z_0\in\Julia(f,U)$, in particular
for repelling and parabolic ones.

By passing to a suitable iterate of $f$, the above definitions are naturally extended to
periodic points.

\subsection{Main results}

Consider a sequence~$\{f_n:U\to\Complex\}_{n\in\Natural}$ of analytic functions with a common
attracting fixed point $z_0\in U$. Suppose that $f_n$ converges to $f_0$ uniformly on each
compact subset of $U$. It is easily follows from arguments of~\cite{Douady} that
$\Basr(z_0,f_n,U)\to\Basr(z_0,f_0,U)$ as to the kernel with respect to~$z_0$ provided $z_0$ is
an attracting or parabolic fixed point of the limit function~$f_0$. At the same time
$\Basr(z_0,f_n,U)$ fails to converge to $\Basr(z_0,f_0,U)$ in general if $z_0$ is a Siegel
point of $f_0$ (see Example~\ref{e_one} in Section~\ref{examples}). Similarly, the dependence
of Julia sets on the function under iteration fails to be continuous at $f_0$ (with respect to
Hausdorff metric) if $f_0$ has (generally speaking, periodic) Siegel points.  Nevertheless, in
the paper~\cite{Kriete} devoted to the continuity of Julia sets for one-parametric families of
transcendental entire functions H.\,Kriete established an assertion, which can be stated as
follows.

\begin{result}\label{KR}
Suppose $f:\Complex\times\Complex\to\Complex$; $(\lambda,z)\mapsto f_\lambda(z)$ is an
analytic family of entire functions $f_\lambda(z)=\lambda z+a_2(\lambda)z^2+\ldots$ and
$\lambda_0:=e^{2\pi i\alpha_0}$, where $\alpha_0\in\Real\setminus\Rational$ is a Diophantine
number. Let $\Delta$ be any Stolz angle at the point $\lambda_0$ with respect to the unit disk
$\{\lambda:|\lambda|<1\}$. Then
$\Basr(0,f_\lambda,\Complex)\to\Basr(0,f_{\lambda_0},\Complex)$ as to the kernel with respect
to~$z_0$ when $\lambda\to\lambda_0$, $\lambda\in\Delta$.
\end{result}

\begin{remark}
It was proved by C.\,Siegel~\cite{Siegel} that for a fixed point with multiplier~$e^{2\pi
i\alpha}$, $\alpha\in\Real\setminus\Rational$, to be a Siegel point, it is sufficient that
$\alpha$ be Diophantine. This condition is not necessary even if restricted to the case of
quadratic polynomials~$f(z):=z^2+c$, $c\in\Complex$ (see \cite[Th.~6]{Brjuno} and
\cite{Yoccoz}). Furthermore, it is easy to construct a nonlinear analytic germ with a Siegel
point for any given~$\alpha\in\Real\setminus\Rational$.
\end{remark}
The Diophantine condition on $\alpha_0$ is substantially employed in~\cite{Kriete}, and in
view of the above remark it is interesting to find out whether this condition is really
essential in Theorem~\ref{KR}. Another question to consider is the role of analytic dependence
of~$f_\lambda$ on~$\lambda$. A possible answer is the following statement improving
Theorem~\ref{KR}.

\begin{theorem}\label{thrm1}
Let $f_0:U\to\Complex$ be an analytic function with a Siegel point ${z_0\in U}$ and
$\{f_n:U\to\Complex\}_{n\in\Natural}$ a sequence of analytic functions with an attracting
fixed point at $z_0$. Suppose that $f_n$ converges to $f_0$ uniformly on each compact subset
of~$U$ and the following conditions hold
\begin{itemize}
\item[(i)] $\big|\arg\big(1-f_n'(z_0)/f_0'(z_0)\big)\big|<\Theta$ for some $\Theta<\pi/2$ and all
$n\in\Natural$;

\item[(ii)] the functions $\big(f_n(z)-f_0(z)\big)/\big(f'_n(z_0)-f'_0(z_0)\big)$,
 $n\in\Natural$, are
uniformly bounded on each compact subset of $U$.
\end{itemize}
Then $\Basr(z_0,f_n,U)$ converges  to $\Basr(z_0,f_0,U)$ as to the kernel with respect to the
point~$z_0$.
\end{theorem}

Condition~(i) in this theorem requires that $\lambda_n:=f_n'(z_0)$ tends to
$\lambda_0:=f'_0(z_0)$ within a Stolz angle, condition~(ii) appears instead of analytic
dependence of~$f_\lambda$ on~$\lambda$, and the Diophantine condition on $\alpha_0$ turns out
to be unnecessary. Both conditions (i) and (ii) are essential. We discuss this in
Section~\ref{examples}.

Dynamics of iterates in the  immediate basin of a fixed point can be described by means of
so-called K\oe nigs function.

Let $z_0$ be a fixed point of an analytic function $f$. The {\it K\oe nigs function}~$\varphi$
associated with the pair $(z_0,f)$ is a solution to the {\it Schr\"oder functional equation}
\begin{equation}\label{eq_Sch}
\varphi\big(f(z)\big)=\lambda \varphi(z),\quad \lambda:=f'(z_0),
\end{equation}
analytic in a neighbourhood of~$z_0$ and subject to the normalization~$\varphi'(z_0)=1$.

It is known~(see e.\,g.~\cite[p.73--76, 116]{Milnor}, \cite{Bargmann}) that the K\oe nigs
function exists, is unique, and can be analytically continued all over~$\Basr(z_0,f,U)$
provided $z_0$ is a geometrically attracting or Siegel fixed point.  If the K\oe nings
function is known, then the iterates can be determined by means of the equality%
\[\varphi\big(f^n(z)\big)=\lambda^n \varphi(z),\quad \lambda:=f'(z_0).\]

By $\varphi_k$, $k\in\Neotr:=\Natural\cup\{0\}$, denote the K\oe nigs function associated with
the pair~$(z_0,f_k)$. We prove the following
\begin{theorem}\label{thrm2}
Under the conditions of Theorem~\ref{thrm1}, the sequence $\varphi_n$ converges to $\varphi_0$
uniformly on each compact subset of~$\Basr(z_0,f_0,U)$.
\end{theorem}
The assertion of Theorem~\ref{thrm2} should be understood in connection with
Theorem~\ref{thrm1}, because the uniform convergence of $\varphi_n$  on a compact set
$K\subset\Basr(z_0,f_0,U)$ requires that $K$ were in the range of definition of $\varphi_n$,
i.e. in $\Basr(z_0,f_n,U)$, for all $n\in\Natural$ apart from a finite number.

\vskip3mm\noindent{\bf Assumption}. Hereinafter it is convenient to assume  without loss of
generality that~$z_0=0$, saving symbol $z_0$ for other purposes.

For any $a\in\Complex$ and $A\subset\Complex$  let us use $aA$ as the short variant of
$\{az:z\in A\}$. By~$\krug{\xi_0}\rho$ denote the disk $\{\xi:|\xi-\xi_0|<\rho\}$, but reserve
the notation~$\UD$ for the unit disk $\krug01$.

\begin{remark}\label{rm_Siegel}
The K\oe nigs function $\varphi_0$ associated with the Siegel point of~$f_0$ admits another
description (see e.g.~\cite[p.116]{Milnor}, \cite{Bargmann}) as the conformal mapping of the
Siegel disk $\Basr(0,f_0,U)$ onto a Euclidean disk $D(0,r)$ that satisfies the condition
$\varphi_0(0)=\varphi'_0(0)-1=0$. From this viewpoint it will be convenient to consider the
conformal mapping $\varphi$, $\varphi(0)=0$, $\varphi'(0)>0$, of $\Basr(0,f_0,U)$ onto the
unit disk~$\UD$ instead of the K\oe nings function $\varphi_0$. Obviously,
$\varphi(z)/\varphi_0(z)$ is constant, and consequently, $\varphi$ satisfies the Schr\"oder
equation~\eqref{eq_Sch} for $f:=f_0$. For shortness, $S$ will stand for $\Basr(0,f_0,U)$. By
$\psi$ denote the inverse function to $\varphi$ and let $S_r:=\psi(r\UD)$,
$\Linija_r:=\partial S_r$ for $r\in[0,1]$. One of the consequences of the fact mentioned above
is that $f_0$ is a conformal automorphism of $S$ and~${S_r,~r\in(0,1)}$.
\end{remark}

During the preparation of this paper another proof of Theorems~\ref{thrm1}\,and\,\ref{thrm2} given in~\cite[p.\,3]{Buff} became known to the author. However, our method allows us also to establish an
asymptotic estimate for the rate of covering level-lines of the Siegel disk by basins of
attraction for one-parametric analytic families. Let $f:W\times U\to\Complex$;
$(\lambda,z)\mapsto f_\lambda(z)$, where $U\ni0$ and $W$ are domains in $\Complex$, be a
family of functions and $\alpha_0$ an irrational number satisfying the following conditions:
\begin{itemize}\label{f_cond}
\item[(i)] $f_\lambda(z)$ depends analytically on both the variable~$z\in U$ and
the parameter~$\lambda\in W$;
\item[(ii)] $f_\lambda(0)=0$ and $f_\lambda'(0)=\lambda$ for all $\lambda\in W$;
\item[(iii)] $\lambda_0:=\exp(2\pi i\alpha_0)\in W$ and the function $f_{\lambda_0}$ has a
Siegel point at $z_0=0$, with $S:=\Basr(0,f_{\lambda_0},U)$ lying in $U$ along with its
boundary $\partial S$.
\end{itemize}

Consider the continued faction expansion of~$\alpha_0$ and denote the $n$-th convergent by
$p_n/q_n$. (See e.g.\cite{Buchstab,Springer} for a detailed exposition on continued factions.) 
For $x>0$ we set
$$n_0(x):=\min\left\{n\in\Natural:\frac{2q_nq_{n+1}}{q_n+q_{n+1}}\ge x\right\},\quad
\ell(x):=q_{n_0(x)}.$$ Notation $\varphi$, $\psi$, $S$, and $S_r$ will refer to the limit function $f_{\lambda_0}$. Lemma~\ref{lm_main} with a slight modification can be used to prove the
following statement.

\begin{theorem}\label{thrmiii}
For any Stolz angle $\Delta$ at the point~$\lambda_0$ there exist a constant~${C>0}$ and a
function~$\varepsilon:(0,1)\to(0,+\infty)$ such that for any~$r\in(0,1)$ the following
statements are true:
\begin{itemize}
\item[(i)] $S_r\subset\Basr(0,f_\lambda,U)$ for all $\lambda\in W\cap\Delta$ satisfying
$|\lambda-\lambda_0|<\varepsilon(r)$;

\item[(ii)] $\varepsilon(r)\ge{C(1-r)^3}/{\ell\big((1-r)^{-\gamma}\big)}$,
\end{itemize}
where $\gamma>\gamma_0:=1+\max\big\{\beta_\psi(1),\,\beta_\psi(-1)\big\}$ and
$\beta_\psi$~stands for the integral means spectrum of the function~$\psi$,
\begin{equation}
\beta_\psi(t):=\limsup_{r\to1-}\frac{\displaystyle\log\int_0^{2\pi}|
\psi'(re^{i\theta})|^t\,d\theta}{-\log(1-r)}.
\end{equation}
\end{theorem}
It is known~\cite{HedenmalmShimorin} that $\beta_\psi(1)\le 0.46$  and $\beta_\psi(-1)\le
0.403$ for any function~$\psi$ bounded and univalent in $\UD$. Consequently, $\gamma_0\le
1.46$.

Theorem~\ref{thrmiii} has been published in~\cite{Sar}. We sketch its proof and specify the function $\varepsilon(r)$ explicitly in Section~\ref{sketch}.

\section{Proof of Theorems}
\subsection{Lemmas}

Denote $\lambda_k:=f'_k(0)$, $k\in\Neotr$. Let us fix arbitrary $n_*\in\Natural$ and consider
the linear family%
\begin{equation}\label{this_family}
f_\lambda[n_*](z):=(1-t)f_0(z)+tf_{n_*}(z),~~
t:=\frac{\lambda-\lambda_0}{\lambda_{n_*}-\lambda_0},\quad z\in U,~\lambda\in\Complex.
\end{equation}%
The number $n_*$ will be not varied throughout the discussion in the present section. So we
shall not indicate dependence on $n_*$ until it is necessary. In particular we shall often
write $f_\lambda$ instead of $f_\lambda[n_*]$.

We need the following elementary statement on approximation of integrals by quadrature sums
(see e.\,g.~\cite[p.\,55--62]{Sobol}).
\begin{result}\label{r_quad}
Suppose $\phi$ is a continuously differentiable function on~$[0,1]$.\linebreak Then for any
$N\in\Natural$ and any set of points $x_0,x_1,\ldots,x_{N-1}\in[0,1]$ the following inequality
holds
\begin{equation}\label{eq_quad}
\left|\int\limits_0^1\phi(x)\,dx-\frac1{N}\sum_{n=0}^{N-1}\phi(x_n)\right|<
Q\big(x_0,x_1,\ldots,x_{N-1}\big)\int\limits_0^1|\phi'(x)|\,dx,
\end{equation}
where $Q\big(x_0,x_1,\ldots,x_{N-1}\big):=
\sup\limits_{x\in[0,1]}\big|F(x;x_0,x_1,\ldots,x_{N-1})-x\big|$ and
\[F(x;x_0,x_1,\ldots,x_{N-1}):=\frac1N\sum_{n=0}^{N-1}\theta(x-x_n),~~~\theta(y):=\left\{
\begin{array}{ll}
1,& \text{\it if~ }y>0,\\%
0,& \text{\it if~ }y\le0.
\end{array}
\right.\]
\end{result}

\begin{remark}\label{nt_QN}
Consider the sequence $x^\beta_n:=\lfr \alpha_0n+\beta\rfr$, where $\lfr\cdot\rfr$ stands for
fractional part, $\alpha_0$ is given by $\lambda_0=e^{2\pi i\alpha_0}$, and $\beta$ is an
arbitrary real number. Denote
\[Q_{\beta,N}:=Q\big(x^\beta_0,x^\beta_1,\ldots,x^\beta_{N-1}\big).\]%
Since $\alpha_0\in\Irrational$, we have (see, e.\,g. \cite[p.\,102--108]{Sobol}) $Q_{\beta,
N}\to0$ as $N\to+\infty$.
\end{remark}

Fix any $r_0\in(0,1)$. The following lemma allows us to determine $\varepsilon_*>0$ such that
$S_{r_0}\subset\Basr(0,f_\lambda,U)$ whenever $|\arg(1-\lambda/\lambda_0)|<\Theta$ and
$|\lambda-\lambda_0|<\varepsilon_*$. In order to state this assertion we need to introduce
some notation.

Denote \[k_0(z):=\frac{z}{(1-z)^2},~z\in\UD,\quad
k_\gamma(z):=e^{i\gamma}k_0(e^{-i\gamma}z),~\gamma\in\Real, \]%
\[u(z):=\frac{f_{n_*}(z)-f_0(z)}{\lambda_{n_*}-\lambda_0},\quad
H(\xi):=1+\frac{\xi\psi''(\xi)}{\psi'(\xi)}, \]%
\[J(t):=\frac {\xi u'(\psi(\xi))\psi'(\xi)-u(\psi(\xi))H(\lambda_0\xi)}
{\lambda_0\xi\psi'(\lambda_0\xi)},~~\xi:=r_0e^{2\pi it}.
\]%

For $\tau\in(0,-\log r_0)$ and  $N\in\Natural$ we put
\[Q_N:=\inf_{\beta\in\Real}Q_{\beta,N},\quad a_N:=2\pi Q_N \int_0^1\left|J(t)\right|\,dt,\]
\[\Lambda_N(\tau,\varepsilon):= \frac{\sqrt{1+2b^2\cos2\vartheta+b^4}-1+b^2}{2
b\cos\vartheta},\quad \varepsilon>0,\] where $\vartheta:=\Theta+\arcsin a_N$,
$b:=\pi\varepsilon N(1-a_N)/(4\tau)$,
\[\varepsilon_N(\tau):=\frac{1-k_\pi(r_*)/k_\pi(r^*)}{\sup\limits_{z\in
S_{r_*}}\left|1-f_{n_*}(z)/f_0(z)\right|}\left|\lambda_{n_*}-\lambda_0\right|,\quad
r_*:=r_0e^{\tau(1-1/N)},~~ r^*:=r_0e^\tau.\]

\begin{lemma}\label{lm_main}
Let $N\in\Natural$ and $\tau\in(0,-\log r_0)$. If $a_N<\sin(\pi/2-\Theta)$, then
$f^N_{\lambda}\big(S_{r_0}\big)\subset S_{r_0}$ for all $\lambda$ such that
$\big|\arg(1-\lambda/\lambda_0)\big|<\Theta$ and $|\lambda-\lambda_0|<\varepsilon_*$, where
$\varepsilon_*:=\varepsilon_N(\tau)\, \Lambda_N\big(\tau,\varepsilon_N(\tau)\big)$.
\end{lemma}

\begin{remark}
In view of Montel's criterion the inclusion $f^N_{\lambda}\big(S_{r_0}\big)\subset S_{r_0}$ in Lemma~\ref{lm_main} implies that $S_{r_0}\subset\Basr(0,f_{\lambda},U)$. We will use this simple fact without reference.
\end{remark}

Lemma~\ref{lm_main} in a slightly different form has been proved in~\cite{Sar}.
We state its proof here for completeness of the discussion. The scheme of the proof is the
following. The main idea is to fix arbitrary~$z_0\in\Linija_{r_0}$ and consider the function
$s_N(\lambda)=s_N(z_0,\lambda):=\varphi(f^N_\lambda(z_0))$. The first step (Lemma~\ref{lm_koltso}) is to determine a neighbourhood
of $\lambda_0$ where $s_N$ is well-defined, analytic and takes values from a prescribed domain
of the form~$\{\xi:\rho_1<|\xi|<\rho_2\}$. The next step (Lemma~\ref{lm_logder}) is to
calculate the value of $(\partial/\partial\lambda) \log s_N(\lambda)$ at~$\lambda=\lambda_0$,
which turns out to be equal to%
\[A_N(z_0):=\frac{s_N'(\lambda_0)}{s_N(\lambda_0)}=\sum_{k=0}^{N-1}
G\big(\lambda_0^k\varphi(z_0)\big),\]%
where $G$ is an analytic function in $\UD$. The concluding step is to use the equality
$\int_0^1G\big(e^{2\pi it}\varphi(z_0)\big)\,dt=1/\lambda_0$ and Theorem~\ref{r_quad} in order
to estimate $|A_N(z_0)|$ and $|\arg A_N(z_0)|$. This allows us to employ a consequence of the
Schwarz lemma (Proposition~\ref{pr_v}) for proving that $|s_N(\lambda)|\le|\varphi(z_0)|$ for any
$\lambda$ satisfying $|\arg(1-\lambda/\lambda_0)|<\Theta$ and
$|\lambda-\lambda_0|<\varepsilon_*$. Since $z_0\in\Linija_{r_0}$ is arbitrary, this
means that $f^N_\lambda\big(S_{r_0}\big)\subset S_{r_0}$ for all such values of $\lambda$.

\begin{lemma}\label{lm_koltso}
Under the conditions of Lemma~\ref{lm_main}, ${s_N(z,\lambda):=\varphi\big(f^N_\lambda(z)\big)}$ is a  well-defined and analytic function for all
$z\in\overline{S_{r_0}}$ and $\lambda\in\krug{\lambda_0}{\varepsilon_N(\tau)}$. Moreover, the
following inequality holds
\begin{equation}\label{eq_koltso}
r_0e^{-\tau}<|s_N(z,\lambda)|<r_0e^\tau,\quad
z\in\Linija_{r_0},~~\lambda\in\krug{\lambda_0}{\varepsilon_N(\tau)}.
\end{equation}
\end{lemma}
\begin{proof}
Let us show that for any $r_1\in(0,1)$, $r_2\in(r_1,1)$ the following inclusion holds
\begin{equation}\label{vkl}
B(z_0,r_1,r_2):=\left\{z:|z-z_0|<|z_0|\big(1-k_\pi(r_1)/k_\pi(r_2)\big)\right\} \subset
S_{r_2}\setminus\overline{S_{r_3}},
\end{equation}
where $z_0\in\Linija_{r_1}$ and $r_3:=r_1^2/r_2$. To this end we remark that for any
$z_0\in\Linija_{r_1}$ the domain $S_{r_2}\setminus\overline{S_{r_3}}$ contains all points~$z$
such that
\begin{equation}\label{ussl}
\big|\log(z/z_0)\big|<\log\big(k_\pi(r_2)/k_\pi(r_1)\big)
\end{equation}
for some of the branches of $\log$. To make sure this statement is true it is sufficient to
employ the following estimate, see e.\,g. \cite[p.\,117, inequal.\,(18)]{Goluzin},
\begin{equation}\label{log_der_est}
\left|\log \frac{z\psi'(z)}{\psi(z)}\right|\le\log\frac{1+|z|}{1-|z|},\qquad z\in\UD,
\end{equation} Owing to~\eqref{log_der_est}, for any rectifiable curve
$\Gamma\subset\overline{S_{r_2}}\setminus S_{r_3}$ that joins $z_0$ with $\Linija_{r_2}$
or~$\Linija_{r_3}$ we have
\begin{multline*}
\int_\Gamma\left|\frac{d
z}{z}\right|=\int_{\varphi(\Gamma)}\left|\frac{\psi'(\xi)}{\psi(\xi)}\right|\,|d\xi|\ge
\int_{\varphi(\Gamma)}\left|\frac{\psi'(\xi)}{\psi(\xi)}\right|\,d|\xi|\\
\ge\min\left\{\int\limits_{r_1}^{r_2}\frac{(1-r)dr}{(1+r)\,r},
\int\limits_{r_3}^{r_1}\frac{(1-r)dr}{(1+r)\,r}\right\}=\log \big(k_\pi(r_2)/k_\pi(r_1)\big).
\end{multline*}

Using the inequality $|\log (1+\xi)|\le -\log(1-|\xi|)$, $\xi\in\UD,$ we conclude that for any
$z\in B(z_0,r_1,r_2)$,
\begin{multline*}
\big|\log\big(z/z_0\big)\big|=\big|\log\big(1+(z-z_0)/z_0\big)\big|
\\\le-\log\big(1-|z-z_0|/|z_0|\big)<\log\big(k_\pi(r_2)/k_\pi(r_1)\big),
\end{multline*}
i.e. all $z\in B(z_0,r_1,r_2)$ satisfy condition~\eqref{ussl}. Therefore inclusion~\eqref{vkl}
holds.

Let $r\in(0,e^{-\tau/N})$. Set $r':=r e^{\tau/N}$ and $r'':=r e^{-\tau/N}$.  Consider an
arbitrary function~$h$ subject to the following conditions: $h$ is analytic in $S$, $h(0)=0$,
and $|h(z)-z|<|z|\big(1-k_\pi(r)/k_\pi(r')\big)$ for all $z\in\overline{S_r}\setminus\{0\}$.

Set $r_1:=|z_0|$, $r_2:=|z_0|e^{\tau/N}$ for some $z_0\in\overline{S_r}\setminus\{0\}$. Since
$k_\pi(x)/k_\pi(xe^{\tau/N})$ increases with $x\in(0,r]$, the Schwarz lemma can be applied to
the function $h(z)-z$ to conclude that $h(z_0)\in B(z_0,r_1,r_2)$ for all $z_0\in
\overline{S_r}\setminus\{0\}$. Therefore~\eqref{vkl} implies the following inclusions
\begin{equation}\label{ash1}
h(\overline{S_r})\subset S_{r'},
\end{equation}
\begin{equation}\label{ash2}
h(\Linija_r)\subset S_{r'}\backslash\overline{S_{r''}}.
\end{equation}

By considering the function $\big(h(z)-z\big)/z$ with $f_{\lambda_0}(w)$ substituted for~$z$
it is easy to check that since the function $f_{\lambda_0}$ is an automorphism of~$S_r$ for
any ${r\in(0,1]}$ (see Remark~\ref{rm_Siegel}), the above argument can be applied to
${h(z):=f_\lambda(f^{-1}_{\lambda_0}(z))}$ for all
${\lambda\in\krug{\lambda_0}{\varepsilon_N(\tau)}}$ and $r\in(0,r_*]$. Thus
\eqref{ash1},\,\eqref{ash2} imply that for any
${\lambda\in\krug{\lambda_0}{\varepsilon_N(\tau)}}$,
\begin{equation}\label{ash11}
f_\lambda\left(\overline{S_{r_j}}\right)\subset S_{r_{j+1}},\qquad j=0,1,\ldots,N-1,
\end{equation}
\begin{equation}\label{ash21}
f_\lambda\left(\overline{S_{r_j}}\backslash S_{r_{-j}}\right)\subset S_{r_{j+1}}\backslash
\overline{S_{r_{-(j+1)}}},\qquad j=0,1,\ldots,N-1,
\end{equation}
where $r_j:=r_0 e^{j\tau/N}$, $j=0,\pm1,\ldots,\pm N$. Applying~\eqref{ash11} repeatedly, we
see that $f^N_\lambda(\overline{S_{r_0}})\subset S_{r_N}$. Similarly,~\eqref{ash21} implies
that $f^N_\lambda(\Linija_{r_0})\subset  S_{r_{N}}\backslash \overline{S_{r_{-N}}}$. The
former means that the function $s_N(z,\lambda)$ is well-defined and analytic for all
$z\in\overline{S_{r_0}}$ and $\lambda\in \krug{\lambda_0}{\varepsilon_N(\tau)}$, while the
latter means that inequality~\eqref{eq_koltso} holds for indicated values of $\lambda$. This
completes the proof of Lemma~\ref{lm_koltso}.
\end{proof}

\begin{lemma}\label{lm_logder}
Under the conditions of Lemma~\ref{lm_koltso}, the following equality holds
\begin{equation}\label{eq_logder}
A_N(z_0):=\left.\frac{\partial\log
s_N(z_0,\lambda)}{\partial\lambda}\right|_{\lambda=\lambda_0}\!\!\!=\,\,\,
\sum_{k=0}^{N-1}G\big(\lambda_0^k\varphi(z_0)\big),\quad z_0\in\Linija_{r_0},
\end{equation}
where \[G(\xi):=\frac{u\big(\psi(\xi)\big)}{\lambda_0\xi\psi'(\lambda_0\xi)}.\]
\end{lemma}
\begin{proof}
Consider the following function of $n+1$ independent variables
\[g_n(z;\lambda_1,\ldots,\lambda_n):=\left\{\begin{array}{cr}\big(f_{\lambda_n}\circ\ldots\circ
f_{\lambda_1}\big)(z),& n\in\Natural,\\z,& n=0.\end{array}\right.\]%
Note that
\begin{align*}A_N(z_0)=\frac{\varphi'\big(f^N_{\lambda_0}(z_0)\big)}
{s_N(z_0,\lambda_0)}\cdot\left.\frac{\partial
g_N(z_0;\lambda,\ldots,\lambda)}{\partial\lambda}\right|_{\lambda=\lambda_0}&&\text{and }\\
\left.\frac{\partial
g_N(z_0;\lambda,\ldots,\lambda)}{\partial\lambda}\right|_{\lambda=\lambda_0}
\!\!\!=\,\,\,\sum_{k=0}^{N-1}g'_{N,k+1}(z_0;\lambda_0,\ldots,\lambda_0),
\end{align*}
where $g'_{n,j}$ stands for $(\partial/\partial\lambda_j)g_n$. Using the equality
\[g_N(z;\lambda_1,\ldots,\lambda_n)=
g_{N-j}\Big(f_{\lambda_j}\big(g_{j-1}(z;\lambda_1,\ldots,\lambda_{j-1})\big);\lambda_{j+1},\ldots,
\lambda_N\Big),\] we get
\[g'_{N,k+1}(z_0;\lambda_0,\ldots,\lambda_0)=\big(f^{N-k-1}_{\lambda_0}\big)'
\big(f^{k+1}_{\lambda_0}(z_0)\big)\cdot u\big(f^k_{\lambda_0}(z_0)\big).\] Schr\"oder
equation~\eqref{eq_Sch} for $f:=f_{\lambda_0}$ allows us to express $f^j_{\lambda_0}$ and
$\big(f^j_{\lambda_0}\big)'$ in terms of $\varphi$ and~$\psi$. Denoting
$z_j:=f_{\lambda_0}^j(z_0)$, $j\in\Neotr$, we
obtain%
\begin{multline*}
g'_{N,k+1}(z_0;\lambda_0,\ldots,\lambda_0)=
\lambda_0^{N-k-1}\psi'\big(\lambda_0^{N-k-1}\,\varphi(z_{k+1})\big) \varphi'(z_{k+1})\,
u(z_k)\\=%
\lambda_0^{N-k-1}\psi'\big(\lambda_0^{N-k-1}\,\varphi(z_{k+1})\big)\,
\frac{u(z_k)}{\psi'\big(\varphi(z_{k+1})\big)}\\
=\lambda_0^{N-k-1}\psi'\big(\lambda_0^N\varphi(z_0)\big)\,
\frac{u\big(\psi\big(\lambda_0^k\varphi(z_0)\big)\big)}{\psi'\big(\lambda_0^{k+1}\varphi(z_0)\big)}.
\end{multline*}

In the same way, we get%
\[\frac{\varphi'\big(f^N_{\lambda_0}(z_0)\big)}
{s_N(z_0,\lambda_0)}=%
\frac{1}{\psi'\big(\lambda_0^N\varphi(z_0)\big)\lambda_0^N\varphi(z_0)}.\]

Now one can combine the obtained equalities to deduce~\eqref{eq_logder}, which completes the
proof.
\end{proof}

\begin{proposition}\label{pr_v}
Let $\tau>0$ and $\Theta\in(0,\pi/2)$. If a function $v(\varsigma)$ is analytic in $\UD$ and
satisfies the following inequalities
\begin{equation}\label{eq_v}
|v(0)|e^{-\tau}<|v(\varsigma)|<|v(0)|e^{\tau},\quad \varsigma\in\UD,
\end{equation}
\begin{equation*}
\vartheta:=\big|\arg \{v'(0)/v(0)\}\big|+\Theta<\pi/2,
\end{equation*}
then the modulus of~~$t:=\pi v'(0)/(4\tau v(0))$ does not exceed\,\, $1$  and  the following
inequality holds
\begin{equation}\label{maj_vyv}|v(\varsigma)|\ge|v(0)|,\quad \varsigma\in\Xi(\rho_0),
\end{equation}
where $\Xi(\rho)$ stands for the circular sector $\{\varsigma:|\mIm\varsigma|\le |\varsigma|
\sin\Theta\le \rho\sin\Theta\}$ and $\rho_0:=\sqrt{\gamma^2+1}-\gamma$,
$\gamma:=(1-|t|^2)/(2|t|\cos\vartheta)$.
\end{proposition}
\begin{proof}
Replacing $v(\varsigma)$ with $v(\varsigma)/v(0)$, we can suppose that $v(0)=1$. The
multi-valued function
\[\phi(\xi):=h\left(\exp\left(\frac{i\pi\log\xi}{2\tau}\right)\right),\quad
h(z):=-i\,\frac{z-1}{z+1},\]%
maps the annulus $\{\xi:e^{-\tau}<|\xi|<e^\tau\}$ conformally onto~$\UD$ (in the sense
of~\cite[p.\,248]{Goluzin}) and satisfies the conditions $\phi(1)=0$, $\phi'(1)>0$. Since the
composition $f:=\phi\circ v$ can be continued analytically along every path in~$\UD$, it
defines an analytic function $f:\UD\to\UD$, $f(0)=0$. By the Schwarz lemma, $|f'(0)|\le1$. Since ${f'(0)=\phi'(1)v'(0)=\pi v'(0)/(4\tau)=t}$,  the first part of Proposition~\ref{pr_v} is proved. To prove the remaining part we note that~\eqref{maj_vyv} is equivalent
to the inequality $\mRe f(\varsigma)\ge0$. Applying the invariant form of the Schwarz lemma to
$f(z)/z$, we obtain
\begin{equation*}\label{maj}
\left|\frac{f(\varsigma)-f'(0)\varsigma}{\varsigma-\overline{f'(0)}f(\varsigma)}\right|
\le|\varsigma|,\quad \varsigma\in\UD,
\end{equation*}
It follows that $f(\varsigma)$ lies in the closed disk of radius
\[R:=\frac{|\varsigma|^2(1-|t|^2)}{1-|t\varsigma|^2}\]%
centred at%
\[\sigma_0:=\frac{t\varsigma(1-|\varsigma|^2)}{1-|t\varsigma|^2}.\]%
Consequently for the inequality $\mRe f(\varsigma)\ge0$ to be satisfied, it is sufficient that
$\mRe\sigma_0\ge R$. An easy calculation leads to the following condition
\begin{equation*}
\cos(\arg t+\arg\varsigma)\ge\frac{|\varsigma|(1-|t|^2)}{|t|(1-|\varsigma|^2)},
\end{equation*}
which is satisfied for all points of the arc%
\[l(\rho):=\{\varsigma:|\mIm\varsigma|\le |\varsigma| \sin\Theta=\rho\sin\Theta\},\quad
\rho\in(0,1),\]%
provided
\begin{equation}\label{maj2}
\cos\vartheta\ge \frac{\rho(1-|t|^2)}{|t|(1-\rho^2)}.
\end{equation}

The right-hand of~\eqref{maj2} increases with $\rho\in(0,1)$ and $\rho:=\rho_0$
satisfies~~\eqref{maj2}. Therefore inequality~\eqref{maj_vyv} holds for all
$\varsigma\in\bigcup_{\rho\in[0,\rho_0]}l(\rho)=\Xi(\rho_0)$. This completes the proof of
Proposition~\ref{pr_v}.
\end{proof}

\begin{proof}\hskip-.2em{\it of Lemma~\ref{lm_main}}\hskip.2em{}
Consider the function~$s_N(z,\lambda)$ introduced in Lemma~\ref{lm_koltso}. This lemma states
that $s_N(z,\lambda)$ is well-defined and analytic for all $z\in \overline{S_{r_0}}$ and
$\lambda\in\krug{\lambda_0}{\varepsilon_{N}(\tau)}$ and satisfies
inequality~\eqref{eq_koltso}. According to Remark~\ref{rm_Siegel},
$f_{\lambda_0}\big(\Linija_{r}\big)=\Linija_{r}$ for all $r\in[0,1)$. Consequently
$|s_N(z,\lambda_0)|=|\varphi(z)|$, $z\in S$. Therefore for any $z_0\in\Linija_{r_0}$ the
function $v(\varsigma):=1/s_N\big(z_0,\lambda_0(1-\varepsilon_N(\tau)\varsigma)\big)$ is
analytic in~$\UD$ and satisfies inequality~\eqref{eq_v}.

Let us employ now Proposition~\ref{pr_v}. To this end we compute the logarithmic derivative of
$v(\varsigma)$ at $\varsigma=0$. By Lemma~\ref{lm_logder},
\[\frac{v'(0)}{v(0)}=\lambda_0\varepsilon_N(\tau)A_N(z_0)=
\lambda_0\varepsilon_N(\tau)\,\sum_{k=0}^{N-1} G\big(\lambda_0^k\varphi(z_0)\big).\]

Consider the sum $E_N:=\sum_{k=0}^{N-1} G(\lambda_0^{k}\varphi(z_0))/N$. It can be regarded as
an approximate value of the integral
\begin{equation*}
E_*:=\int\limits_0^1G\big(r_0e^{2\pi i(t+t_0)}\big)\,dt,
\end{equation*}
where $t_0\in\Real$ is an arbitrary number, which does not affect $E_*$:
\[%
E_*=\frac{1}{2\pi i}\!\!\int\limits_{|\xi|=r_0}\!\!\!\frac{G(\xi)}{\xi}\,d\xi=
\res\limits_{\xi=0}\frac{G(\xi)}\xi=G(0)=\frac1{\lambda_0}.
\]%
Applying Theorem~\ref{r_quad} to the points $x_n:=x_n^\beta$,
 $\beta:=(\arg\varphi(z_0))/(2\pi)-t_0$, and the function $\phi(t):=G\big(r_0e^{2\pi
i(t+t_0)}\big)$, we get the following estimate%
\[|E_N-E_*|< Q_{\beta,N}\int_0^1\big|(d/dt) G\big(r_0e^{2\pi i(t+t_0)}\big)\big|\,dt.\] Since
$t_0\in\Real$ is arbitrary real, we have
\begin{equation}\label{ots0}%
|E_N-E_*|\le Q_{N}\int_0^1\big|(d/dt) G\big(r_0e^{2\pi i(t+t_0)}\big)\big|\,dt.
\end{equation}%
The function under the sign $\int_0^1|\cdot|dt$ is
\[%
\frac{d G\big(r_0e^{2\pi i(t+t_0)}\big)}{dt}=2\pi i\, \xi G'(\xi)=2\pi i\,J(t+t_0),\quad
\xi:=r_0e^{2\pi i(t+t_0)}.
\]%

From~\eqref{ots0} it follows that
\begin{equation*}
\left|\frac1N\,\cdot\sum_{k=0}^{N-1}
G(\lambda_0^{k}\varphi(z_0))-\frac{1}{\lambda_0}\right|\le a_N,
\end{equation*}
and hence,
\begin{equation}\label{l_eq}
\left|\frac1{N}\cdot\frac{v'(0)}{v(0)}-\varepsilon_{N}(\tau)\right|\le a_N
\varepsilon_{N}(\tau).
\end{equation}
Since by condition $0\le a_N<1$, inequality~\eqref{l_eq} implies that
\begin{equation}\label{gar1}
\left|\frac{v'(0)}{v(0)}\right|\ge N(1-a_N)\varepsilon_{N}(\tau),
\end{equation}
\begin{equation}\label{gar2}
\left|\arg\frac{v'(0)}{v(0)}\right|\le\arcsin a_N.
\end{equation}

Now if we recall that validity of~\eqref{eq_v} has been already verified and take into
account~\eqref{gar1},\,\eqref{gar2}, we see that the conditions of Proposition~\ref{pr_v} are
satisfied. Therefore, by elementary reasoning we see that~\eqref{maj_vyv} holds for
all~$\varsigma\in\Xi\big(\Lambda_N\big(\tau,\varepsilon_N(\tau)\big)\big)$. In terms of $s_N$
this means that
\begin{equation}\label{sq}
|s_N(z_0,\lambda)|\le|s_N(z_0,\lambda_0)|=r_0,\quad \lambda\in\Xi_0,
\end{equation}
where \[\Xi_0:=\big\{\lambda:\big|\lambda-\lambda_0\big|<\varepsilon_*,
\big|\arg(1-\lambda/\lambda_0)\big|<\Theta\big\}.\]

Since $z_0\in\Linija_{r_0}=\partial S_{r_0}$ is arbitrary in the above arguments, by the
maximum modulus theorem, inequality~\eqref{sq} implies that $|\varphi(f^N_{\lambda}(z))|<r_0$
for all $z\in S_{r_0}$ and $\lambda\in\Xi_0$. Therefore for indicated values of~$\lambda$ we
have $f^N_\lambda\big(S_{r_0}\big)\subset S_{r_0}$. This completes the proof of
Lemma~\ref{lm_main}.
\end{proof}

\subsection{Proof of Theorem~\ref{thrm1}} Suppose that the sequence $\{f_n\}_{n\in\Natural}$
satisfies the conditions of Theorem~\ref{thrm1}. Then every subsequence of $f_n$ also meets
these conditions. So we have only to prove that $S:=\Basr(0,f_0,U)$ is the kernel of the
sequence $A_n:=\Basr(0,f_0,U)$, that is:
\begin{itemize}
\item[(i)] any compact set $K\subset S$ lies in all but finite number of $A_n$;
\item[(ii)] $S$ is the largest domain that contains the point $z=0$ and satisfies
condition~(i).
\end{itemize}

Now we employ Lemma~\ref{lm_main} in order to prove~(i). To this end we should fix any
$r_0\in(0,1)$ such that $S_{r_0}\supset K$, specify appropriate values of $N$ and $\tau$, and
trace the dependence on the choice of $n_*$. As a result we would prove that
\begin{equation}\label{eq_inf_eps}
\varepsilon^0_*:=\inf_{n_*\in\Natural}\varepsilon_*>0.
\end{equation}
Since $\lambda_n\to\lambda_0$ as $n\to+\infty$, \eqref{eq_inf_eps}  would imply  that 
 ${K\subset S_{r_0}\subset\Basr(z_0,f_n,U)}$ for all $n\in\Natural$ large enough.

Set $\tau:=(1+r_0)/(2r_0)$. In view of condition~(ii) of Theorem~\ref{thrm1},
\[L:=\sup_{n_*\in\Natural}\left(\int_0^1\left|J(t)\right|\,dt\right)<+\infty.\]%
Since by Remark~\ref{nt_QN}, $Q_N\to0$ as $N\to+\infty$, there exists $N\in\Natural$ such that
\begin{equation*}
Q_N<\frac{\sin\left(\pi/4-\Theta/2\right)}{2\pi L}.
\end{equation*}
Fix any such value of $N$. Then $a_N<\sin(\pi/4-\Theta/2)<\sin(\pi/2-\Theta)$. Hence
Lemma~\eqref{lm_main} is applicable to the specified values of $N$ and $\tau$.

Let us estimate $\varepsilon_*$ from below.  In view of condition~(ii) of Theorem~\ref{thrm1},
\[\varepsilon_0:=\inf_{n_*\in\Natural}\varepsilon_N(\tau)>0.\] Denote
$b:=\pi\varepsilon_N(\tau)N(1-a_N)/(4\tau)$, $b_1:=\min\{1,b\}$.\\ Since
$\vartheta=\Theta+\arcsin a_N<\pi/4+\Theta/2<\pi/2$, we have
\begin{multline*}\Lambda_N\big(\tau,\varepsilon_N(\tau)\big)\ge
\frac{\sqrt{1+2b_1^2\cos2\vartheta+b_1^4}-1+b_1^2}{2b_1\cos\vartheta}\\\ge
\frac{\sqrt{1+2b_1^2\cos2\vartheta+(b_1^2\cos2\vartheta)^2}-1+b_1^2}{2b_1\cos\vartheta}
\\>b_1\cos\vartheta> b_1\cos(\pi/4+\Theta/2)\\\ge
\cos(\pi/4+\Theta/2)\min\left\{1,
\frac{\pi\varepsilon_0N\big(1-\sin(\pi/4-\Theta/2)\big)}{4\tau}\right\}=:C_0.
\end{multline*}
The constant $C_0$ is positive and does not depend on $n_*$. From the inequality
${\varepsilon_*>\varepsilon_0C_0}$ it follows that \eqref{eq_inf_eps} takes place. This proves
assertion~(i).

To prove~(ii) let us assume the converse. Then there exists a domain $S'\not\subset S$, $0\in
S'$, satisfying~(i). Let $z_0\in S'\setminus S$ and $\Gamma\subset S'$ be a curve that joins
points $z=0$ and $z_0$. Consider any domain $D$ such that $\Gamma\subset D$ and $K:=\overline
D\subset S'$. By the assumption, $K\subset A_n$ for all $n$ large enough. Now we claim that
\begin{equation}\label{eq_E}
D\subset E(f_0,U).
\end{equation}
Consider an arbitrary $\zeta_0\in D$. Suppose that $\zeta_0\not\in E(f_0,U)$. Then there
exists $j_0\in\Natural$ such that $f_0^{j_0}$ is well-defined (and so analytic) in some domain
$D_0\ni\zeta_0$, $D_0\subset D$, with $f^{j}_0(\zeta_0)\in U$, $j<j_0$, but
$f^{j_0}_0(\zeta_0)\not\in U$. Since the sequence $f_n$ converges to $f_0$ uniformly on each
compact subset of~$U$, the sequence $f_n^{j_0}$ converges to $f_0^{j_0}$ uniformly on each
compact subset of~$D_0$. According to Hurwitz's theorem, this means that
$f_n^{j_0}(D_0)\not\subset U$ for all $n\in\Natural$ large enough. Consequently, $D\not\subset
E(f_n,U)$ for large~$n$. At the same time, $K=\overline D\subset A_n\subset E(f_n,U)$ for all~$n$ large enough. This contradiction proves~\eqref{eq_E}.

The remaining part of the proof  depends on the properties of the domain~$U$. Since
$U\subset\Complex$, we have three possibilities:
\begin{itemize}
\item[(Hyp)] The domain~$U$ is hyperbolic. Then by Montel's criterion, $\Fatou(f_0,U)$ coincides
with the interior of $E(f_0,U)$. Since $D\ni 0$ is connected, we conclude that
$z_0\in\Gamma\subset D\subset S$. With this fact contradicting the assumption, the proof
of~(ii) for the hyperbolic case is completed.

\item[(Euc)] The domain~$U$ coincides with $\Complex$. The functions $f_n$, $n\in\Neotr$, are
entire functions.

\item[(Cyl)] The domain~$U$ is the complex plane punctured at one  point.
\end{itemize}

Let us prove (ii) for case~(Euc). Since $\Gamma\cap\partial S\neq\emptyset$ and $\partial
S\subset\Julia(f_0,\Complex)$, we have $D\cap\Julia(f_0,\Complex)\neq\emptyset$.
The classical result proved for entire functions by I.N.\,Baker~\cite{BakerR} asserts that the
Julia set coincides with the closure of the set of all repelling periodic points. Therefore,
$D$ contains a periodic point of $f_0$ different from~$0$. Owing to Hurwitz's theorem, the
same is true for $f_n$ provided $n$ is large enough. This leads to a contradiction, because
the immediate basin of attraction $\Basr(0,f_n,U)$ contains no periodic points except for the
fixed point~$z=0$. Assertion~(ii) is now proved for case~(Euc).

It remains to consider case~(Cyl). Similarly to case~(Euc), we need only to show that
$D\setminus\{0\}$ contains a periodic point. By means of linear transformations we can assume
that $U=\Complex\setminus\{1\}$. From~\eqref{eq_E} it follows that functions
\[\phi_n(z):=\frac{f_0^n(z)-z}{f_0^n(z)-1},\quad n\in\Natural,\]%
does~not assume values $1$ and $\infty$ in $D$. Since $D\cap\Julia(f_0,U)\neq\emptyset$, the
family $\{\varphi_n\}_{n\in\Natural}$ is not normal in $D$. Hence, due to Montel's criterion,
there exists $z_1\in D$ and $n_0\in\Natural$ such that $\phi_{n_0}(z_1)=0$ and so $z_1\in D$
is a periodic point of~$f_0$. This completes the proof of~(ii) for case~(Cyl).

By now  (i) and (ii) are shown to be true. Theorem~\ref{thrm1} is  proved.\qed

\subsection{Proof of Theorem~\ref{thrm2}} Fix any $r_0\in(0,1)$. As in the proof of
Theorem~\ref{thrm1}  one can make use of Lemma~\ref{lm_main} to show that there exist $n_1,
N\in\Natural$ such that $f^N_n(S_{r_0})\subset S_{r_0}$ for all~$n>n_1$. By
Remark~\ref{rm_Siegel} the function~$\varphi_0$ maps $S$ conformally onto a Euclidian disk
centred at the origin. It is convenient to rescale the dynamic variable, by replacing
$f_k$, $k\in\Neotr$, with $r f_k(z/r)$ for some constant~$r>0$, so that $\varphi_0(S)=\UD$ (or
equivalently $\varphi_0=\varphi$). Then the functions ${g_n(\zeta):=(1/r_0)\big(\varphi_0\circ
f_n^N\circ\varphi_0^{-1}\big)(r_0\zeta)}$, $n>n_1$, are defined and analytic in~$\UD$.
Furthermore, $g_n(0)=0$ and $g_n(\UD)\subset\UD$ for any $n>n_1$. Let us observe that for any
analytic function $f$ with a geometrically attractive or Siegel fixed point~$z_0$ the K\oe
nigs function~$\varphi$ associated with the pair~$(z_0,f)$ is the same as that of the
pair~$(z_0,f^N)$. Hence it is easy to see that the function
$\phi_n(\zeta):=\varphi_n\big(\varphi_0^{-1}(r_0\zeta)\big)/r_0$ is the K\oe nigs function
associated with~$(0,g_n)$. Since $S=\varphi_0^{-1}(\UD)$ and $r_0\in(0,1)$ is arbitrary, it
suffices to prove that~$\phi_n(\zeta)\to\zeta$ as $n\to+\infty$ uniformly on each compact
subset of~$\UD$.

According to Remark~\ref{rm_Siegel}, the function $f_0$ is a conformal automorphism of~$S$.
Therefore, with $f_n$ converging to $f_0$ uniformly on each compact subset of $U\supset S$,
there exists $n_2\ge n_1$ such that for all $n>n_2$ functions $f^N_n$ and consequently $g_n$
are univalent in $S_{r_0}$ and in~$\UD$, respectively. It follows~(see e.g.\cite{Gor}) that
$\phi_n$, $n>n_2$, are also univalent in~$\UD$. The convergence of $f_n$ to $f_0$ implies also
that $g_n$ converges to~$g_0$, $g_0(\zeta):=\lambda_0^N\zeta$, uniformly on each compact
subset of~$\UD$.

We claim that there exists a sequence $\big\{r_n\in(0,1)\big\}_{n\in\Natural}$ converging
to~$1$ such that for all $n>n_2$ the domain $\phi_n(r_n\UD)$ is contained in some
disk~${\{\xi:|\xi|<R_n\}}$ that lies in $\phi_n(\UD)$. Owing to the Carath\'eodory convergence
theorem and normality of the family~$\{\phi_n:n\in\Natural,n>n_2\}$, this statement would
imply convergence of the sequence $\phi_n$ to the identity map and hence the proof of
Theorem~\ref{thrm2} would be completed.

By $p/q$ and $p'/q'$ let us denote some successive convergents of the number\linebreak
${\alpha_n:=\big(\arg g'_n(0)\big)/(2\pi)=(\arg\lambda_n^N)/(2\pi)}$ (regardless of whether
$\alpha_n$ is irrational or not). Put $\Omega_n:=\phi_n(\UD)$, $\kappa_n:=-\log|g_n'(0)|=-N\log|\lambda_n|$,
$a_n:=\kappa_n(q-1)$, and $b_n:=\pi(1/q+2/q')$. Consider a point
$\zeta_0\in\UD$ and make use of the following inequality (see e.\,g. \cite[p.\,117,
inequal.\,(18)]{Goluzin}) from the theory of univalent function
\[\left|\log\frac{\zeta\phi_n'(\zeta)}{\phi_n(\zeta)}
\right|\le\log\frac{1+|\zeta|}{1-|\zeta|}, \quad\zeta\in\UD,\]%
to obtain
\begin{equation}\label{eq_Gamma}
\int_\Gamma\left|\frac{dw}{w}\right|\ge-\log\big(4k_\pi(|\zeta_0|)\big),\quad
k_\pi(z):=\frac{z}{(1+z)^2},~z\in\UD,
\end{equation}
where $\Gamma$ is any rectifiable curve that joins $\xi_0:=\phi_n(\zeta_0)$
with~$\partial\Omega_n$ and lies in~$\Omega_n$ except for one of the endpoints. The equality
in~\eqref{eq_Gamma} can occur only if $\phi_n$ is a rotation of the Koebe function
$k_0(z):=z/(1-z)^2$ and $\Gamma$ is a segment of a radial half-line. It follows that
$\Omega_n$ contains the annular sector%
\[\Sigma:=\big\{\xi_0e^{x+iy}:|x|\le a_n,\, |x|\le
b_n,\,x,y\in\Real\big\}\]%
provided $|\zeta_0|\le r_n:=k_\pi^{-1}\big((1/4)\exp(-\,\sqrt{a_n^2+b_n^2})\big).$ Moreover,
$\Omega_n$ is invariant under the map $\zeta\mapsto\lambda_n^N\zeta$. Indeed,
\[\lambda_n^N\zeta=\lambda_n^N\phi_n\big(\phi_n^{-1}(\zeta)\big)=
\phi_n\big(g_n\big(\phi_n^{-1}(\zeta)\big)\big)\in\Omega_n\]%
for all $\zeta\in\Omega_n.$ Denote 
$$\Sigma_0:=\big\{\xi_0e^{x+iy}:{|x|\le \pi/q,}\, |x|\le b_n,\,x,y\in\Real\big\},\quad \lambda_*:=e^{-\kappa+2\pi ip/q}$$
Since $p$ and $q$ are coprime integers, the union of the annular
sectors $\lambda_*^j\Sigma_0$, $j=0,1,\ldots,q-1$, contains the circle~$\xi_0\UC$, $\UC:=\partial\UD$. The inequality from the theory of continued fractions $|\alpha_n-p/q|\le1/(qq')$ implies that
\[\lambda_*^j\Sigma_0\subset\big(\lambda_n^{N}\big)^j\Sigma,\quad j=0,1,\ldots,q-1.\]

Therefore, for any $\xi_0\in\phi_n(r_n\UD)$ the domain~$\Omega_n$ contains the
circle~$\xi_0\UC$. It follows that $\phi_n(r_n\UD)$ is a subset of some
disk~$\{\xi:|\xi|<R_n\}$ contained in~$\Omega_n$.

It remains to choose the successive convergents $p/q$ and $p'/q'$ of $\alpha_n$ in such a way
that $r_n\to1$ as $n\to+\infty$. To this end we fix some successive convergents $p/q$ and
$p'/q'$ of $\alpha_*:=(\arg\lambda_0^N)/(2\pi)$ and note that $p/q$ and $p'/q'$ are also
successive convergents of $\alpha_n$ provided $n$ is large enough, because
$\alpha_n\to\alpha_*$ as $n\to+\infty$. Using the fact that $\kappa_n\to0$ as $n\to+\infty$
 and that the denominators of convergents of the irrational number
$\alpha_*$ forms unbounded increasing sequence, we see that it is possible to choose $p/q$ for
each $n$ in such a way that $\sqrt{a^2_n+b^2_n}\to0$ and, consequently, $r_n\to1$  as
$n\to+\infty$.

The proof of Theorem~\ref{thrm2} is now completed. \qed

\section{Proof of Theorem~\ref{thrmiii}}\label{sketch}
In this section we sketch the proof of Theorem~\ref{thrmiii}. First of all we note that the proof of Lemma~\ref{lm_main} does not use the fact that the dependence of~$f_\lambda[n_*]$ (see equation~\eqref{this_family}) on the parameter $\lambda$ is linear. So Lemma~\ref{lm_main} can be applied to any analytic family~$f_\lambda$ satisfying conditions~(i)\,--\,(iii) on page~\pageref{f_cond}, provided some notations are modified to a new (more general) setting. First of all we have to redefine $$u(z):=\left.\frac{\partial f_\lambda(z)}{\partial\lambda}\right|_{\lambda=\lambda_0}$$ Then fix any $r\in(0,1)$ and consider the modulus of continuity of the family ${h_\lambda:=f_\lambda/f_{\lambda_0}}$ calculated at $\lambda=\lambda_0$,
$$\omega_r(\delta):=\sup\Big\{\big|1-f_{\lambda}(z)/f_{\lambda_0}(z)\big|:z\in S_r,\lambda\in
W\cap \krug{\lambda_0}\delta\Big\},\quad\delta>0.$$ This quantity, as a function of $\delta$, is defined, continuous, and increasing on the interval~$I^*:=(0,\delta^*)$, $\delta^*:=\dist(\lambda_0,\partial W)$, with $\lim_{\delta\to+0}\omega_r(\delta)=0$. Therefore there exists an inverse function~$\omega^{-1}_r:(0,\epsilon^*)\to(0,+\infty)$, where
$\epsilon^*:=\lim_{\delta\to\delta^*-0}\omega_r(\delta)$. If $\epsilon^*\neq+\infty$, then we set
$\omega_r^{-1}(\epsilon):=\delta^*$ for all $\epsilon\ge\epsilon^*$. Now we can redefine $\varepsilon_N(\tau)$ as
$$\varepsilon_N(\tau):=\omega^{-1}_{r_*}\big(1-k_\pi(r_*)/k_\pi(r^*)\big),\quad
r_*:=r_0e^{\tau(1-1/N)},~~ r^*:=r_0e^\tau.$$
Finally, define $\Theta$ to be equal to the half-angle of~$\Delta$. To apply Lemma~\ref{lm_main} we need the following
\begin{proposition}\label{Q_prop}
For any $n\in\Natural$ the following inequality holds
\begin{equation}\label{Fbeta}
Q_{q_n}<(1/q_n+1/q_{n+1})/2.
\end{equation}
\end{proposition}
\begin{proof}
Fix $n\in\Natural$. Due to the inequality $|\alpha_0-p_n/q_n|<1/(q_n q_{n+1})$
there exists $\gamma\in(0,1/q_{n+1})$ such that \begin{equation}\label{Buch_eq}|\alpha_0-p_n/q_n|<\gamma/q_n.\end{equation}
Let $\beta_0:=\big(1/q-(-1)^n\gamma\big)/2$. Taking into account that $p_n$ and $q_n$ are coprime integers one can deduce by means of the inequalities $\gamma<1/q_{n+1}<1/q_n$, $(-1)^n(\alpha_0-p_n/q_n)>0$, and~\eqref{Buch_eq} that 
\begin{equation}\label{beta0}
Q_{\beta_0,q_n}<(1/q_n+1/q_{n+1})/2.
\end{equation}
This proves the proposition.
\end{proof}

Now let us show how Theorem~\ref{thrmiii} can be proved. Fix $r_0\in(0,1)$. Define~$\varepsilon(r_0)$ in the followin way.

According to Proposition~\ref{Q_prop}, $0<a_N<\sin\big(\pi/4-\theta/2\big)$ for
\begin{gather*}\label{setNtau}
N:=\ell\left(\left.{2\pi}{\textstyle\int_0^1\big|J(t)\big|\,dt}\right/
{\sin\big(\pi/4-\theta/2\big)}\right),\quad \tau:=\log\frac{1+2r_0}{3r_0}.
\end{gather*}
Hence, Lemma~\ref{lm_main} can be used with the specified values of $N$ and $\tau$. Therefore, we can set $\varepsilon(r_0):=\varepsilon_*$, so that statement~(i) in Theorem~\ref{thrmiii} becomes true. Let us show that statement~(ii) of this theorem is also true, assuming that $r_0$ is sufficiently close to~$1$.

Since $\overline S\subset U$ there exists $C_1>0$ such that $$\left|1-\frac{f_\lambda(\psi(\xi))}{f_{\lambda_0}(\psi(\xi))}\right|
<C_1|\lambda-\lambda_0|$$ for all $\xi\in\UD$ and $\lambda\in\krug{\lambda_0}{\varepsilon^0}$, where $\varepsilon^0>0$ is choisen so that $\overline{\krug{\lambda_0}{\varepsilon^0}}\subset W$. It follows that
\begin{equation*}
\omega^{-1}_r(s)\ge\min\big\{\varepsilon^0,s/C_1\big\},\qquad s>0,~r\in(0,1).
\end{equation*}
Elementary calculations show that
\begin{equation*}
1-\frac{k_\pi(r_*)}{k_\pi(r^*)}\ge 1-\exp\left(-\frac{\,\tau(1-r^*)}{N(1+r^*)}\right)
\ge C_2\frac{(1-r_0)^2}{N}
\end{equation*}
for some constant $C_2>0$. Combining these two inequalities we obtain
\begin{equation*}
\varepsilon_N(\tau)\ge C_3\frac{(1-r_0)^2}{N},\qquad C_3:=C_2/C_1.
\end{equation*}
We estimate $\Lambda_N(\tau,\varepsilon_N(\tau))$ in the same way as in the proof of Theorem~\ref{thrm1} to conclude that
\begin{equation*}
\varepsilon_*\ge C\,\frac{(1-r_0)^3}{N}
\end{equation*}
for some constant $C>0$. To complete the proof we use the following inequalities (see e.\,g.~\cite[p.\,52]{Goluzin}):
\begin{gather*}
\label{H}
\left|\frac{\xi\psi''(\xi)}{\psi'(\xi)}-\frac{2r^2}{1-r^2}\right|\le
\frac{4r}{1-r^2},\quad 0\le r=|\xi|<1,\\
\label{Der}
\left|\frac{\psi'(\xi)}{\psi'(0)}\right|\ge\frac{1-r}{(1+r)^3},~~~\qquad\qquad
0\le r=|\xi|<1,
\end{gather*}
which imply that $N\le\ell\big((1-r_0)^{-\gamma}\big)$ for all $r_0<1$ sufficiently close to~$1$.
\qed

\section{Essentiality of conditions in Theorem~\ref{thrm1}}\label{examples}
In this section we show that conditions~(i) and~(ii) in Theorem~\ref{thrm1} are essential. As
for condition~(i) this can be regarded as a consequence of lower semi-continuity of the Julia
set.
\begin{example}\label{e_one}
Consider the family $f_\lambda(z):=\lambda z+z^2$ in the whole complex plane ($U:=\Complex$).
The map $\lambda\mapsto\Julia(f_\lambda,\Complex)$ is lower semi-continuous~\cite{Douady},
i.e.
\[\Julia(f_{\lambda_*},\Complex)~\subset~%
\bigcap\limits_{\varepsilon>0}~~~\bigcup\limits_{\delta>0}~~
\bigcap\limits_{|\lambda-\lambda_*|<\delta}~O_\varepsilon\big(\Julia(f_\lambda,\Complex)%
\big)~~~ \text{ for any } \lambda_*\in\Complex,\]%
where $O_\varepsilon(\cdot)$ stands for the $\varepsilon$-neighbourhood of a set. Let
$\lambda_0:=e^{2\pi i\alpha_0}$, $\alpha_0\in\Irrational$, and $\alpha_n\in\Rational$ converge
to $\alpha_0$ as $n\to+\infty$. The point~$z_0:=0$ is a parabolic fixed point
of~$f_{\lambda_n^0}$, $\lambda_n^0:=\exp(2\pi i\alpha_n)$, and so
$0\in\Julia(f_{\lambda_n^0},\Complex)$. Due to lower semi-continuity
of~$\lambda\mapsto\Julia(f_\lambda,\Complex)$ at the points $\lambda^0_n$, there exists a
sequence ${\{\mu_n\in(0,1)\}_{n\in\Natural}}$ such that
$\krug{0}{1/n}\cap\Julia(f_{\lambda_n},\Complex)\neq\emptyset$, $\lambda_n:=\mu_n\lambda_n^0$,
$n\in\Natural$. It follows that $\Basr(0,f_{\lambda_n},\Complex)\to\{0\}$ as to the kernel.
Assume that $f_{\lambda_0}$, $\lambda_0:=\exp(2\pi i\alpha_0)$, has a Siegel point at $z_0=0$.
This is  the case if $\alpha_0$ is a Brjuno number (\hskip.05em\cite[Th.~6]{Brjuno}, see
also~\cite{Yoccoz}). The sequence $f_n:=f_{\lambda_n}$ satisfies all  conditions of
Theorem~\ref{thrm1} except for condition~(i), but the conclusion of Theorem~\ref{thrm1} fails
to be true. Therefore condition~(i) is an essential one.
\end{example}

It is known~\cite[p.\,44]{BuffHab} that condition~(ii) can be omitted in Theorem~\ref{thrm1}
provided that the multiplier of the Siegel fixed point $\lambda_0:=f_0'(z_0)$ equals to
$\exp(2\pi i \alpha_0)$ for some Brjuno number~$\alpha_0$. However, if no such assumptions
concerning~$\alpha_0$ are made, condition~(ii) cannot be omitted. This fact is demonstrated by
the following example

\begin{example}\label{e_two}
 Let $\alpha_0$ be an irrational real number. By $q_n$ denote the denominator of the
$n$-th convergent of~$\alpha_0$. Consider the sequence of polynomials
\[f_n(z):=\frac{\lambda_0\left(z+z^{q_n+1}\right)}{\hbox{$1+1/2^{q_n}$}},\quad
\lambda_0:=e^{2\pi i\alpha_0},\] converging to $f_0(z)=\lambda_0 z$ uniformly on each compact
subset of~$\UD$.

We claim that the sequence of domains $\Basr(0,f_n,\UD)$ does not converge to ${\Basr(0,f_0,\UD)=\UD}$ as to the kernel, provided the growth rate of $q_n$ is sufficiently high. Assume the converse. Then for all $n\in\Natural$ large enough, say for $n>n_0$,
the inclusion $D_{48}\subset\Basr(0,f_n,\UD)$ holds, where $D_j:=j/(j+1)\UD$, $j\in\Natural$.
It follows that $f_n^m(D_{48})\subset\UD$ for all $n>n_0$, $n\in\Natural$, and $m\in\Natural$.
Hence the family $\Phi:=\{f_n^m\}_{n>n_0,\,n,m\in\Natural_0}$ is normal in the disk~$D_{48}$.
In particular, there exist constants $C_1>1$, $C_2>0$ such that
\begin{gather}\label{ofp}
\left|\left(f^m_n\right)'(z)\right|<C_1, \quad z\in D_8,\quad n>n_0,~n,m\in\Natural_0,
\\\label{osp}\left|\left(f^m_n\right)''(z)\right|<C_2,
\quad z\in D_8,\,\quad n>n_0,~ n,m\in\Natural_0.
\end{gather}
Furthermore, by Schwarz lemma,
\begin{equation}\label{Schw}
f^m_n(D_4)\subset D_5,~~~f^m_n(D_6)\subset D_7 \quad n>n_0,~n,m\in\Natural_0,
\end{equation}

Consider functions $\displaystyle g_n:=f_n^{q_n}$, $\tilde g_n:=\tilde f_n^{q_n}$,%
\[\tilde
f_n(z):=\frac{\exp\left(2\pi i
p_n/q_n\right)}{\hbox{$1+1/2^{q_n}$}}\left(z+z^{q_n+1}\right),\quad z\in\UD, \quad n>n_0,\quad
n\in\Natural,\]%
where $p_n$ stands for the numerator of the $n$-th convergent of~$\alpha_0$. Apply the
following inequality
\begin{multline}\label{ff}
\big|\tilde f_n(z)-f_n(z)\big|=\big|f_n(z)\big|\cdot\big|\lambda_0-\exp\left(2\pi i
p_n/q_n\right)\big|\\\le
4\pi\left|\alpha-\frac{p_n}{q_n}\right|\le\frac{4\pi}{q_nq_{n+1}},\quad z\in\UD,
\end{multline}
to prove that
\begin{equation}\label{gg}
\big|\tilde g_n(z)-g_n(z)\big|<\frac{4\pi C_1}{q_{n+1}},\quad z\in D_4,
\end{equation}
for all $n\in\Natural$ large enough.

Since $q_n\to+\infty$ as $n\to+\infty$, there exists $n_1\in\Natural$, $n_1\ge n_0$, such that
\[\frac{4\pi}{q_nq_{n+1}}<\frac1{72}~~\text{and}~~\frac{4\pi C_1}{q_{n+1}}<\frac1{42},~\quad
n>n_1,\quad n\in\Natural.\]%
We shall show that for all $n>n_1$, $n\in\Natural$, and $k=1,2,\ldots,q_n-1$ the following
implication holds
\begin{equation}\label{imp}
\Big(P(j),~~j=1,2,\ldots,k\Big)\Longrightarrow P(k+1),
\end{equation}
\newbox\mybox
\setbox\mybox=\hbox{\hskip4mm\parbox{11cm}{ \it  $\big|\tilde f_n^{j-1}(z)\big|<1$, $z\in
D_4$, and}}
\newbox\fixbox
\setbox\fixbox=\vbox{where $P(j)$: \it $\big|\tilde f_n^{j-1}(z)\big|<1$, $z\in D_3$, and}
\newdimen\fixdim
\fixdim=\ht\mybox\advance\fixdim by-\ht\fixbox%

\noindent where $P(j):=$ \Big[\lower\fixdim\box\mybox
\begin{equation}\label{req}
\big|\tilde f_n^{j}(z)-f_n^{j}(z)\big|<\frac{4j\pi C_1}{q_nq_{n+1}},\quad z\in D_4.~~~\Big]
\end{equation}

Now let $n>n_1$ and $P(j)$ take place for all $j=1,2,\ldots,k$. Relations~\eqref{Schw},
\eqref{ff}, and~\eqref{req} imply the following inclusions
\begin{equation}\label{vlozh}
\tilde f_n(D_6)\subset D_8,\qquad \tilde f_n^j(D_4)\subset D_6, \quad j=1,2,\ldots,k.
\end{equation}
For $j:=k$ the latter guarantees that $|\tilde f_n^k(z)|<1$, $z\in D_4$. Fix any~$z\in D_4$
and denote $w_j:=\tilde f_n^{j}(z)$, $\tilde \xi_j:=\tilde f_n(w_j)$, $\xi_j:=f_n(w_j)$.
According to~\eqref{Schw}~and~\eqref{vlozh}, we have $w_j\in D_6$, $\tilde \xi_j, \xi_j\in
D_8$, $j=1,2,\ldots,k$. Taking this into account, from~\eqref{ofp}~and~\eqref{ff}, we get the
following inequality
\begin{multline*}
\Big|\tilde f_n^{k+1}(z)-f_n^{k+1}(z)\Big|\le%
\sum_{j=0}^{k}\Big|\big(f_n^{k-j}\circ \tilde f_n^{j+1}\,\big)(z)- \big(f_n^{k-j+1}\circ
\tilde
f_n^j\,\big)(z)\Big|\\=%
\sum_{j=0}^{k}\Big|\big(f_n^{k-j}\circ \tilde f_n\,\big)(w_j)- \big(f_n^{k-j}\circ
f_n\,\big)(w_j)\Big|\\=%
\sum_{j=0}^{k}\Big|f_n^{k-j}(\tilde\xi_j)- f_n^{k-j}(\xi_j)\Big|\\<%
\sum_{j=0}^kC_1|\tilde \xi_j-\xi_j|\le%
\frac{4(k+1)\pi C_1}{q_n q_{n+1}}.
\end{multline*}
Therefore, \eqref{req} holds also for $j:=k+1$. This proves implication~\eqref{imp}.

For $j:=1$ inequality~\eqref{req} follows from~\eqref{ff}. Hence $P(1)$ is valid. Owing
to~\eqref{imp}, $P(1)$ implies $P(q_n)$. Therefore, inequality~\eqref{gg} holds for all
$n>n_1$.

The functions $\tilde g_n$  have the fixed point $\tilde z_*:=1/2$. Now we apply~\eqref{gg} to
show that if
\begin{equation}\label{qqn}
q_{n+1}\ge 2^{q_n},\qquad n\in\Natural,
\end{equation}
then for any sufficiently large $n\in\Natural$ the function $g_n$ has also a fixed point
$z_*\in D_3\setminus\{0\}$. Straightforward calculation gives
\begin{equation*}\label{eq_pr1}
\tilde g_n'(\tilde
z_*)=l_n:=\left(\frac{1+(q_n+1)/2^{q_n}}{\hbox{$1+1/2^{q_n}$}}\right)^{q_n}>1.
\end{equation*}

From~\eqref{osp},\,\eqref{gg}, and the Cauchy integral formula it follows that
\begin{equation*}\label{osp_tilde}
\left|\tilde g_n''(z)\right|<C_3:=C_2+51200\pi C_1/q_{n+1}, \quad z\in D_3,\quad n>n_1.
\end{equation*}   Now we assume that $n\in\Natural$ is large enough and
apply Rouch\'e's theorem to the functions $\tilde g_n(z)-z$ and $g_n(z)-z$ in the disk
$B_n:=\big\{z:|z-1/2|<\rho_n\big\}$, where $\rho_n:=(l_n-1)/(2C_3)$.  Since $B_n\subset D_3$,
we have $$\mRe
 \frac d{dz}\big(\tilde g_n(z)-z\big)>\frac{l_n-1}2\,,\qquad z\in B_n.$$ It follows that
\begin{equation}\label{circ}
|\tilde g_n(z)-z|\ge \frac{(l_n-1)\rho_n}2,\qquad z\in \partial B_n.
\end{equation}

Inequalities~\eqref{gg},\,\eqref{qqn}, and \eqref{circ} imply that $|\tilde g_n(z)-z|>|\tilde
g_n(z)-g_n(z)|$ for all $z\in\partial B_n$. Consequently, $g_n(z)-z$ vanishes at some point
$z_*\in B_n$. At the same time, the immediate basin $\Basr(0,f_n,\UD)$ contains no periodic
points of $f_n$ except for the fixed point at $z_0=0$. Therefore,
$D_3\not\subset\Basr(0,f_n,\UD)$ for large $n$. This fact implies that the sequence
$\Basr(0,f_n,\UD)$ does not converge to~$\UD$ as to the kernel.

It easy to see that the prescribed sequence~$f_n$ satisfies all  conditions of
Theorem~\ref{thrm1} with $U:=\UD$ except for condition~(ii), but the conclusion fails to hold.
This shows that (ii) is also an essential condition in Theorem~\ref{thrm1}.
\end{example}

{

\end{bibliography}}
\small
\vskip15mm
\noindent{\sc Department of Mathematics, University of Bergen, Johannes\\ Brunsgate 12, 5008 Bergen, NORWAY}

\vskip1mm \noindent {\sc Department of Mechanics and Mathematics, Saratov State\\
University,
  Astrakhanskaya 83, 410012 Saratov, RUSSIA}\\[2mm]
{\it E-mail address:} \href{mailto:Pavel.Gumenyuk@math.uib.no}{\tt Pavel.Gumenyuk@math.uib.no}, \href{mailto:gumenuk@sgu.ru}{\tt gumenuk@sgu.ru}

\end{document}